\documentclass[12pt,a4paper]{article}

\usepackage[title]{appendix}
\usepackage{a4wide}
\usepackage{epsfig}
\usepackage{authblk}
\usepackage{amsmath}
\usepackage{amssymb}
\usepackage{graphicx}
\usepackage{amsthm,dsfont,amsfonts,fancyhdr}

\usepackage[running]{lineno}

\usepackage{pst-node}

\usepackage{tikz}
\usepackage{enumerate}

\usepackage[bookmarks=false,hyperfootnotes=false]{hyperref}

\usepackage{comment}

\newcommand{\QR}{\square_q}
\newcommand{\NQR}{\overline{\square}_q}

\newtheorem{propo}{Proposition}[section]
\newtheorem{defi}[propo]{Definition}

\newtheorem{lemma}[propo]{Lemma}

\newtheorem{theo}[propo]{Theorem}
\newtheorem{examp}[propo]{Example}
\newtheorem{remar}[propo]{Remark}
\newtheorem{ques}[propo]{Question}

\newtheorem{cor}[propo]{Corollary}
\newtheorem{prop}[propo]{Proposition}

\newcommand{\bl}{\begin{lemma}}
\newcommand{\el}{\end{lemma}}

\def\Cay{{\rm Cay}}

\usepackage{indentfirst,latexsym,bm}

\def\Sym{{\rm Sym}}

\def\Aut{{\rm Aut}}

\def\K{{\rm K}}

\def\HA{{\rm HA}}
\def\HS{{\rm HS}}
\def\HC{{\rm HC}}
\def\AS{{\rm AS}}
\def\TW{{\rm TW}}
\def\SD{{\rm SD}}
\def\CD{{\rm CD}}
\def\PA{{\rm PA}}

\def\H{{\rm H}}

\def\soc{{\rm soc}}

\begin{document}
\title{Comparing classes of highly symmetric graphs: From $2$-arc-transitive to $2$-distance-transitive}

\author[a,b]{Wei Jin}
\author[c,d]{Jack H. Koolen\textsuperscript{*}}
\author[c]{Chenhui Lv}

\affil[a]{\footnotesize{School of Mathematics and Computational Science, Key Laboratory of Intelligent Computing and Information Processing of Ministry of Education, Xiangtan University, Xiangtan, Hunan, 411105, P.R.China}}

\affil[b]{\footnotesize{School of Statistics and Data Science, 
Jiangxi University of Finance and Economics, Nanchang, Jiangxi, 330013, P.R.China}}

\affil[c]{\footnotesize{School of Mathematical Sciences, University of Science and Technology of China, Hefei, Anhui, 230026, P.R.China}}

\affil[d]{\footnotesize{CAS Wu Wen-Tsun Key Laboratory of Mathematics, University of Science and Technology of China, Hefei, Anhui, 230026, P.R.China}}


\maketitle

\begin{center}
\small
\textsuperscript{*}Corresponding author: Jack H. Koolen
(\texttt{koolen@ustc.edu.cn}).\\
E-mail addresses: \texttt{jinweipei82@163.com} (W. Jin),
\texttt{lch1994@mail.ustc.edu.cn} (C. Lv).
\end{center}

\begin{abstract}
A $2$-distance-transitive graph is a vertex-transitive graph whose vertex stabilizer is transitive on both the first- and second-step neighborhoods. 
Every distance-transitive graph of diameter at least $2$, as well as every non-complete $2$-arc-transitive graph, is $2$-distance-transitive.
In this paper, we first determine the vertex-quasiprimitive types of $2$-distance-transitive graphs of odd order, partially answering a question posed by A. Devillers, M. Giudici, C. H. Li and C. E. Praeger in 2012.
We then prove that a $2$-distance-transitive graph of valency $p+1$, where $p$ is a prime, is $2$-arc-transitive if and only if it has girth at least $4$. We also show that every locally-primitive $2$-distance-transitive graph of valency at most $8$ is $2$-arc-transitive, with the icosahedron as the unique exception. Finally, we prove that if $\Gamma$ is a $G$-locally-primitive, $(G,2)$-distance-transitive graph of valency at least $3$ and $G$ is soluble, then either $\Gamma\cong \K_{p,p}$ for some prime $p$, or the order of $\Gamma$ is not square-free.
\end{abstract}

\vspace{2mm}

\hspace{-17pt}{\bf Keywords:} $2$-distance-transitive graphs, $2$-arc-transitive graphs, locally-primitive graphs.

\hspace{-17pt}\textbf{Mathematics Subject Classification (2020):} 05E18; 20B25

\section{Introduction}
In this paper, all graphs are finite, simple, connected and undirected. Undefined notation and terminology will be introduced in Section~2.

A graph $\Gamma$ is said to be \emph{$s$-distance-transitive} if, for each $i\leq s$, where $s$ is an integer not exceeding the diameter of $\Gamma$, the automorphism group of $\Gamma$ is transitive on the set of ordered pairs of vertices at distance $i$.
The study of $s$-distance-transitive graphs goes back to Higman~\cite{DGH-1}.
Devillers et al.~\cite{DGLP-ldt} studied locally $s$-distance-transitive graphs using the normal quotient strategy developed for $s$-arc-transitive graphs in~\cite{Praeger-4}.
Corr, Schneider and the first author~\cite{CJS-2dt} investigated $2$-distance-transitive but not $2$-arc-transitive graphs of valency at most $5$. Two years later, the family of $2$-distance-transitive Cayley graphs over cyclic groups was determined precisely in~\cite{CJL-2019}.
$2$-distance-transitive graphs with diameter $2$ are strongly regular graphs of rank $3$, which have been well classified (see, e.g., \cite[Section~11]{BM-2022}).
For more work on $2$-distance-transitive graphs, see \cite{DJLP-clique,FengKwak-2p3,JT-2016,Praeger-3,qdk-1}.

A graph $\Gamma$ is called \emph{$2$-arc-transitive} if its automorphism group acts transitively on the sets of arcs and $2$-arcs. 
The first remarkable results on $2$-arc-transitive graphs are due to Tutte~\cite{Tutte-1,Tutte-2}; since then, this family of graphs has been extensively studied; see
\cite{ACMX-1996,CP-1983,GLP-1,IP-1,Li-abeliancay-2008,PLY-2010}.
By definition, every non-complete $2$-arc-transitive graph is $2$-distance-transitive, but the converse does not hold in general.
Indeed, if a $2$-distance-transitive graph has girth $3$, then it is not $2$-arc-transitive; the icosahedron provides such an example.
Thus, the family of non-complete $2$-arc-transitive graphs is properly contained in the family of $2$-distance-transitive graphs.

Lemma~5.3 of~\cite{DGLP-ldt} (see Lemma~\ref{dt-quotient12}) provides a global analysis of the family of locally $(G,s)$-distance-transitive graphs with $s \geq 2$.
By this result, each $(G,2)$-distance-transitive graph other than $\K_{m[b]}$ with $m\geq 3$ and $b\geq 2$ has a connected $(G/N,s)$-distance-transitive quotient graph, where $s=\min\{2,d(\Gamma_N)\}$, corresponding to a normal subgroup $N$ of $G$ such that $G/N$ acts quasiprimitively or bi-quasiprimitively on the vertex set of the quotient graph.
Thus, in order to classify $(G,2)$-distance-transitive graphs, we need to know all the `basic' quotient graphs, including complete graphs, for which the corresponding group acts quasiprimitively or bi-quasiprimitively on the vertex set.

For $(G,s)$-distance-transitive graphs with $s \geq 2$, Devillers, Giudici, Li, and Praeger~\cite{DGLP-ldt} posed the following question.

\begin{ques} [{cf.~\cite[Question 6.6]{DGLP-ldt}}]\label{2dt-basic-question1}
{\rm Let $\Gamma$ be a $(G,s)$-distance-transitive graph of valency at least $3$, where $s\geq 2$.
Suppose that $G$ is quasiprimitive on $V(\Gamma)$.
What quasiprimitive types arise for $G$ on $V(\Gamma)$?
In particular, can $G$ act quasiprimitively on $V(\Gamma)$ with a type that does not arise for
$(G,s)$-arc-transitive graphs?
}
\end{ques}

For $(G,2)$-arc-transitive graphs of odd order and valency at least $3$, Li~\cite{LCH-odd-2001} proved that if $G$ is quasiprimitive on the vertex set, then the quasiprimitive action must be of type \AS.
For $(G,2)$-distance-transitive graphs of odd order where $G$ acts quasiprimitively on the vertex set, our first theorem determines all possible vertex-quasiprimitive action types and hence answers Question~\ref{2dt-basic-question1} in the odd-order case.

\begin{theo}\label{2dt-qp-th1}
Let $\Gamma$ be a $(G,2)$-distance-transitive graph of odd order.
If $G$ is quasiprimitive on $V(\Gamma)$ of type $X$, then $X$ is one of the following three types: \AS, \PA, or \HA. Moreover, there are infinitely many examples for each type.
\end{theo}

A graph $\Gamma$ is called \emph{locally-primitive} if for each vertex $u$, the vertex stabilizer $\Aut(\Gamma)_u$ is primitive on the neighborhood of $u$. 
The study of the family of  well-known locally-primitive graphs is an active topic; see~\cite{LLP-abeliancay-2011,LPM-prime-2009,LPVZ-2002,Potocnik-abeliancay}.
Every non-complete $2$-arc-transitive graph is locally $2$-transitive, and hence both locally-primitive and $2$-distance-transitive.
Therefore, the class of locally-primitive $2$-distance-transitive graphs can be regarded as an intermediate class between non-complete $2$-arc-transitive graphs and $2$-distance-transitive graphs.
For $G$-locally-primitive, $(G,2)$-distance-transitive graphs of odd order, Theorem~\ref{2dt-qp-th1} leads to the following corollary.

\begin{cor}\label{2dt-odd-cor1}
Let $\Gamma$ be a $G$-locally-primitive, $(G,2)$-distance-transitive graph of odd order with valency at least $3$.
If $G$ is quasiprimitive on $V(\Gamma)$ of type $X$, then $X$ is either of type \AS\ or of type \PA.
\end{cor}

\begin{remar}
{\rm Combining Li~\cite{LCH-odd-2001}, Theorem~\ref{2dt-qp-th1}, and Corollary~\ref{2dt-odd-cor1}, we observe a clear hierarchy of admissible types $X$ for graphs of odd order and valency at least $3$ where $G$ acts quasiprimitively on $V(\Gamma)$:
\begin{itemize}
\item[(1)] If $\Gamma$ is $(G,2)$-arc-transitive, then Li~\cite{LCH-odd-2001} showed that $X$ is of type \AS.
\item[(2)] If $\Gamma$ is $G$-locally-primitive and $(G,2)$-distance-transitive, then Corollary~\ref{2dt-odd-cor1} shows that $X$ is either of type \AS\ or of type \PA.
\item[(3)] If $\Gamma$ is $(G,2)$-distance-transitive, then Theorem~\ref{2dt-qp-th1} shows that $X$ is of type \AS, \PA, or \HA.
\end{itemize}
}
\end{remar}

As observed above, there is a natural inclusion relationship among these three classes of graphs. Consequently, it is of great interest to investigate and compare the relative sizes of these classes under additional specific restrictions, such as bounded valency. In Theorem~\ref{2dt-theo-prime1}, by restricting the valency to $p+1$ where $p$ is a prime, we compare the class of $2$-distance-transitive graphs with the class of $2$-arc-transitive graphs. In Theorem~\ref{2dt-localp-th1}, under the assumption that the valency is at most $8$, we compare the class of locally-primitive $2$-distance-transitive graphs with the class of $2$-arc-transitive graphs.

Let $\Gamma$ be a $2$-distance-transitive graph.
If $\Gamma$ has girth $3$, then $\Gamma$ has two types of $2$-arcs, and hence it is not $2$-arc-transitive.
If $\Gamma$ has girth at least $5$, then for each pair of vertices at distance $2$, there is a unique $2$-arc between them, and hence $\Gamma$ is $2$-arc-transitive.
In the case where $\Gamma$ has girth $4$, the graph $\Gamma$ may or may not be $2$-arc-transitive. 

Our second theorem characterizes $(G,2)$-distance-transitive graphs of valency $p+1$ where $p$ is a prime.
For a characterization of $2$-geodesic-transitive graphs of prime valency, see \cite{DJLP-prime}.

\begin{theo}\label{2dt-theo-prime1}
Let $\Gamma$ be a connected $2$-distance-transitive graph of valency $p+1$, where $p$ is a prime. Then $\Gamma$ is $2$-arc-transitive if and only if it has girth at least $4$.
\end{theo}

The third theorem compares the classes of locally-primitive $2$-distance-transitive graphs and $2$-arc-transitive graphs under the assumption that the valency is at most $8$.

\begin{theo}\label{2dt-localp-th1}
Let $\Gamma$ be a locally-primitive $2$-distance-transitive graph of valency at most $8$. 
Then $\Gamma$ is either isomorphic to the icosahedron or $2$-arc-transitive.
\end{theo}

In the study of highly symmetric graphs, the interplay between the order of a graph and the algebraic structure of its automorphism group is a central theme. Graphs of square-free order are particularly interesting, since permutation groups acting on sets of square-free size are highly constrained. This leads to a natural question: can such restrictive arithmetic conditions on the number of vertices, together with strong local symmetry, be compatible with a structurally ``simple'' group, such as a soluble group?
Our final theorem addresses this question. It gives a strict structural constraint: if a graph is $G$-locally-primitive and $(G,2)$-distance-transitive, then a soluble group $G$ is generally too restrictive to support square-free order. Apart from a family of complete bipartite graphs, such highly symmetric graphs of square-free order must arise from insoluble groups.

\begin{theo}\label{2dt-lp-sf}
Let $\Gamma$ be a $G$-locally-primitive, $(G,2)$-distance-transitive graph of valency at least $3$.
If $G$ is soluble, then either $\Gamma\cong \K_{p,p}$ for some prime $p$, or the order of $\Gamma$ is not square-free.
\end{theo}

In the proof of Theorem~\ref{2dt-localp-th1}, we utilize \cite[Corollary 4.8]{JKL-2025-1} to show that every amply regular graph with diameter $4$ and parameters $k = 7$ and $\mu = 3$ has exactly $32$ vertices. Resolving this specific parameter case was actually one of the primary motivations for our recent study in~\cite{JKL-2025-1}. In that paper, we constructed three infinite families of amply regular graphs with diameter $d = 4$ and parameters $(4q+4, q, 0, \frac{q - 1}{2})$, arising from Paley graphs, Peisert graphs, and Paley digraphs, respectively. To further connect our previous combinatorial constructions with the highly symmetric graphs studied in the present paper, we include Appendix~\ref{app-2at}. In this appendix, we prove that two of these families---namely, those constructed from Paley graphs and Paley digraphs---are in fact $2$-arc-transitive.

\textbf{Outline of the paper.} In the next section, we provide preliminaries and definitions.  In Section~3, we study $(G,2)$-distance-transitive graphs of odd order where $G$ acts quasiprimitively, and prove Theorem~\ref{2dt-qp-th1} and Corollary~\ref{2dt-odd-cor1}. Section~4  investigates $2$-distance-transitive graphs of valency $k = p+1$, leading to the proof of Theorem~\ref{2dt-theo-prime1}. 
In Section~5, we first establish restrictions on the parameters of $2$-distance-transitive graphs and then use them to analyze locally-primitive $2$-distance-transitive graphs of valency at most $8$, proving Theorem~\ref{2dt-localp-th1}. Section~6 studies $(G,2)$-distance-transitive graphs of square-free order when $G$ is soluble; by analyzing the Fitting subgroup of $G$, we prove Theorem~\ref{2dt-lp-sf}. Finally, Appendix~\ref{app-2at} shows that the two families of amply regular graphs constructed from Paley graphs and Paley digraphs in~\cite{JKL-2025-1} are $2$-arc-transitive.


\section{Preliminaries}
In this section, we introduce several definitions from group and graph theory that will be used throughout the paper. For group-theoretic terminology not defined here, we refer the reader to \cite{Cameron-1,DM-1,Wielandt-book}.

\subsection{Basic definitions in graph theory}
For a graph $\Gamma$, we denote its \emph{vertex set} by $V(\Gamma)$, its \emph{edge set} by $E(\Gamma)$, and its \emph{automorphism group} by $\Aut(\Gamma)$. The \emph{order} of $\Gamma$ is the number of its vertices, that is, $|V(\Gamma)|$. For two vertices $x, y \in V(\Gamma)$, if there exists an edge between them, we say that $x$ is adjacent to $y$, or that $x$ and $y$ are neighbors; this is denoted by $x \sim y$.

The \emph{distance} $d(x, y) = d_\Gamma(x, y)$ between two vertices $x, y \in V(\Gamma)$ is the length of a shortest path connecting them in $\Gamma$. The \emph{diameter} $d = d(\Gamma)$ of $\Gamma$ is the maximum distance among all pairs of vertices, and the \emph{girth} of $\Gamma$ is the length of a shortest cycle. For each $x \in V(\Gamma)$ and $0 \leq i \leq d$, let $\Gamma_i(x)$ denote the set of vertices at distance $i$ from $x$. Additionally, we define $\Gamma_{-1}(x) = \Gamma_{d+1}(x) = \emptyset$. For convenience, we abbreviate $\Gamma_1(x)$ as $\Gamma(x)$. The \emph{valency} of a vertex $x$ in $\Gamma$ is the cardinality $d(x) = d_\Gamma(x) := |\Gamma(x)|$. A graph $\Gamma$ is called \emph{regular} of valency $k$ if $|\Gamma(x)| = k$ for all $x \in V(\Gamma)$.

Let $\K_m$ denote the \emph{complete graph} on $m$ vertices, and let $t\K_m$ denote the disjoint union of $t$ copies of $\K_m$. The \emph{complete bipartite graph} with parts of sizes $m$ and $n$ is denoted by $\K_{m,n}$.
Let $\K_{m[n]}$ denote the \emph{complete multipartite graph} with $m$ parts, each of size $n$.

The \emph{complement} of a graph $\Gamma$, denoted by $\overline{\Gamma}$, is the graph with the same vertex set as $\Gamma$, where two vertices are adjacent in $\overline{\Gamma}$ if and only if they are not adjacent in $\Gamma$.

The \emph{induced subgraph} of a graph $\Gamma$ on a vertex subset $U \subseteq V(\Gamma)$, denoted by $[U]$, is the graph whose vertex set is $U$ and whose edge set consists of all edges of $\Gamma$ with both endpoints in $U$. The graph $[\Gamma(x)]$ is called the \emph{local graph} of $\Gamma$ at $x$.

The \emph{bipartite double} of a graph $\Gamma$ is defined as the graph with vertex set $V(\Gamma) \times \{0, 1\}$, where two vertices $(x, i)$ and $(y, j)$ are adjacent if and only if $x \sim y$ in $\Gamma$ and $i \neq j$. For simplicity, we write $(x, 0)$ as $x^+$ and $(x, 1)$ as $x^-$.

\subsection{Amply regular graphs and distance-regular graphs}
An \emph{amply regular graph} with parameters $(v, k, \lambda, \mu)$ is a $k$-regular graph on $v$ vertices such that every pair of adjacent vertices has exactly $\lambda$ common neighbors, and every pair of vertices at distance two has exactly $\mu$ common neighbors.

A \emph{$(0,2)$-graph} is a connected graph in which each pair of distinct vertices have either $0$ or $2$ common neighbors.

A \emph{strongly regular graph} with parameters $(v, k, \lambda, \mu)$ is a regular graph of order $v$ and valency $k$ such that every pair of adjacent vertices has exactly $\lambda$ common neighbors, and every pair of non-adjacent vertices has exactly $\mu$ common neighbors. It is equivalent to being an amply regular graph of diameter $2$ with the same parameters. For further properties of strongly regular graphs, see~\cite{BM-2022} and~\cite[Chapter~10]{GR}.

A connected graph $\Gamma$ with diameter $d$ is called \emph{distance-regular} if there exist integers $b_i$ and $c_i$ for $0 \leq i \leq d$ such that, for each pair of vertices $x, y \in V(\Gamma)$ with $d(x, y) = i$, the vertex $y$ has exactly $c_i$ neighbors in $\Gamma_{i-1}(x)$ and $b_i$ neighbors in $\Gamma_{i+1}(x)$, where $b_d = c_0 = 0$ (cf.~\cite[p.126]{BCN}). It is straightforward to verify that $c_1 = 1$ always holds. In particular, a distance-regular graph  is regular of valency $k := b_0$, and for notational convenience, we define $a_i := k - b_i - c_i$. The constants $a_i$, $b_i$, and $c_i$ are called the \emph{intersection numbers}, and the array
$(b_0 = k, b_1, \dots, b_{d-1}; c_1, \dots, c_d)$ is referred to as the \emph{intersection array}.

We present two constraints on the parameters of amply regular graphs that will be used repeatedly.

\begin{theo} [{cf.~\cite[Theorem 1.5.5]{BCN}}]\label{Thm1.5.5}
Let $\Gamma$ be a connected amply regular graph with parameters $(v,k,\lambda,\mu)$
and diameter $d \geq 3$. Let $k_2 = \frac{k(k-\lambda-1)}{\mu}$.
Then
\[
k_2 \geq k \geq \lambda + \mu + 1.
\]
Moreover, $k_2 = k$ (equivalently, $k = \lambda + \mu + 1$), if and only if
$\Gamma$ is a polygon or a Taylor graph.
\end{theo}


\begin{lemma} [{cf.~\cite[Lemma~2.1]{JKL-2025-1}}]\label{b1geqk} 
Let $\Gamma$ be a connected amply regular graph with parameters $(v,k,\lambda,\mu)$
and diameter $d \geq 3$. Let $b_1 = k - \lambda - 1$.
Then $b_1 \geq \frac{1}{3}(k+1)$ holds.
\end{lemma}

\subsection{Paley graphs and Taylor graphs}
Let $q$ be a prime power such that $q \equiv 1 \pmod{4}$. The \emph{Paley graph} $P(q)$ is defined as the graph with vertex set $\mathbb{F}_q$, where two distinct vertices $x, y \in \mathbb{F}_q$ are adjacent if and only if $x - y$ is a nonzero square in $\mathbb{F}_q$. The congruence condition guarantees that $-1$ is a square in $\mathbb{F}_q$, ensuring the graph is undirected. The graph $P(q)$ is a strongly regular graph with parameters $(v, k, \lambda, \mu) = \left(q, \frac{q - 1}{2}, \frac{q - 5}{4}, \frac{q - 1}{4} \right)$.

A distance-regular graph $\Gamma$ with intersection array $(k, \mu, 1;\ 1, \mu, k)$ is called a \emph{Taylor graph}. Such a graph has diameter $3$ and order $2(k + 1)$. Moreover, for each vertex $x \in V(\Gamma)$, we have $|\Gamma(x)| = |\Gamma_2(x)| = k$ and $|\Gamma_3(x)| = 1$.

Given a graph $\Gamma$ with vertex set $V(\Gamma)$, its \emph{Taylor double} is the graph with vertex set $\{x^\varepsilon \mid x \in V(\Gamma),\ \varepsilon = \pm 1\}$, where two distinct vertices $x^\delta$ and $y^\varepsilon$ are adjacent (for $x \neq y$) if and only if $\delta\varepsilon = 1$ when $x \sim y$, and $\delta\varepsilon = -1$ when $x \not\sim y$. Given a strongly regular graph $\Delta$ with $v$ vertices satisfying $k = 2\mu$, its \emph{Taylor extension} is the Taylor double of the graph $\{\infty\} + \Delta$, where $\infty$ is a new vertex adjacent to every vertex of $\Delta$. The resulting graph is a Taylor graph on $2(v + 1)$ vertices with intersection array $\{v, v - k - 1, 1;\ 1, v - k - 1, v\}$ (see~\cite[p.19]{BM-2022}).

\subsection{Permutation groups and  quotient graphs}
Let $G$ be a permutation group acting on a set $\Omega$. Recall that a \emph{$G$-invariant partition} of $\Omega$ is a partition $\mathcal{B}=\{B_1,B_2,\ldots,B_n\}$ such that, for all $g\in G$ and all $i$, we have $B_i^g=B_j$ for some $j$. The parts of such a partition are often called \emph{blocks}. They satisfy the defining property that $B\cap B^g$ is either empty or equal to $B$, for every block $B$ and every $g\in G$. A partition $\mathcal{B}$ of $\Omega$ is said to be \emph{nontrivial} if $2 \leq |\mathcal{B}| < |\Omega|$; otherwise, $\mathcal{B}$ is said to be \emph{trivial}. In particular, if $N\lhd G$, then the set of $N$-orbits in $\Omega$ forms a $G$-invariant partition.

Let $\Gamma$ be a graph and let $G\leq \Aut(\Gamma)$. Suppose that $\mathcal{B} = \{B_1, B_2, \ldots, B_n\}$ is a partition of $V(\Gamma)$. The \emph{quotient graph} $\Gamma_{\mathcal{B}}$ is the graph whose vertex set is $\mathcal{B}$, in which two parts $B_i$ and $B_j$ are adjacent if and only if there exist vertices $x \in B_i$ and $y \in B_j$ such that $x$ is adjacent to $y$ in $\Gamma$.
Whenever $\mathcal{B}$ is the set of $N$-orbits in $V(\Gamma)$ for some $N\lhd G$, we also write $\Gamma_{\mathcal{B}}=\Gamma_N$ and call it a \emph{$G$-normal quotient} of $\Gamma$.
The graph $\Gamma$ is said to be a \emph{cover} of $\Gamma_{\mathcal{B}}$ if, for each pair of adjacent blocks $B_i$ and $B_j$ and each $x \in B_i$, we have
$|\Gamma(x) \cap B_j| = 1$.

A permutation group $G$ acting on a set $\Omega$ is said to be
\emph{$2$-transitive} if it is transitive on the set of ordered pairs of distinct points in $\Omega$.

\subsection{$s$-Distance-transitive and $s$-arc-transitive}
Let $G \leq \Aut(\Gamma)$ and let $s \leq d(\Gamma)$. We say that $\Gamma$ is \emph{locally $(G,s)$-distance-transitive} if, for every vertex $v \in V(\Gamma)$, the stabilizer $G_v$ acts transitively on $\Gamma_i(v)$ for all $i \leq s$; if $s = d(\Gamma)$, then $\Gamma$ is said to be \emph{locally $G$-distance-transitive}. 

A $(G,s)$-distance-transitive graph is a locally $(G,s)$-distance-transitive graph such that $G$ is transitive on $V(\Gamma)$; when $s = d(\Gamma)$, we say that $\Gamma$ is \emph{$G$-distance-transitive}. If $G = \Aut(\Gamma)$, we omit $G$ and simply say that $\Gamma$ is \emph{$s$-distance-transitive}; in particular, when $s = d(\Gamma)$, we say that $\Gamma$ is \emph{distance-transitive}.

By definition, if $s > d(\Gamma)$, then $\Gamma$ is not $(G, s)$-distance-transitive.  
For example, the complete graph $\K_n$ is not $(G, 2)$-distance-transitive for any group $G$.

In the study of $(G, s)$-distance-transitive graphs, the following constants are essential.  
Their definition is inspired by the notion of intersection numbers for distance-regular graphs.

\begin{defi}
{\rm
Let $\Gamma$ be a $(G, s)$-distance-transitive graph, and let $u \in V(\Gamma)$ and $v \in \Gamma_i(u)$ with $0 \leq i \leq s$. Then the number of neighbors of $v$ in $\Gamma_{i-1}(u)$, $\Gamma_i(u)$, and $\Gamma_{i+1}(u)$ is independent of the choice of $v$, and these numbers are denoted by $c_i$, $a_i$, and $b_i$, respectively.}
\end{defi}

It is clear that $a_i + b_i + c_i$ equals the valency of $\Gamma$ whenever these constants are well-defined.  
In particular, for $(G, 2)$-distance-transitive graphs, the constants $c_i$, $a_i$, and $b_i$ are always well-defined for $i = 0, 1, 2$, and hence such graphs are amply regular.

A sequence $(u_0, u_1, \ldots, u_s)$ of vertices in a graph is called an \emph{$s$-arc} if $u_i$ is adjacent to $u_{i+1}$ for all $i \in \{0, \ldots, s-1\}$ and $u_i \neq u_{i+2}$ for all $i \in \{0, \ldots, s-2\}$.
A graph $\Gamma$ is said to be \emph{$(G, s)$-arc-transitive} if $G$ acts transitively on the vertex set and on the set of $s$-arcs of $\Gamma$.  
If $G = \Aut(\Gamma)$, then we simply say that $\Gamma$ is \emph{$s$-arc-transitive}.

A triple $(u, v, w)$ of vertices in a graph $\Gamma$ is called a \emph{$2$-geodesic} if $v$ is adjacent to both $u$ and $w$, and $u \neq w$ and $u$ is not adjacent to $w$.  
A $G$-arc-transitive graph that is not complete is said to be \emph{$(G, 2)$-geodesic-transitive} if $G$ acts transitively on the set of $2$-geodesics.  
If $G = \Aut(\Gamma)$, then we simply say that $\Gamma$ is \emph{$2$-geodesic-transitive}.

By definition, every non-complete $2$-arc-transitive graph is $2$-geodesic-transitive, and every $2$-geodesic-transitive graph is $2$-distance-transitive.  
However, neither implication is reversible in general.

We will need the following result on locally $(G,s)$-distance-transitive graphs from \cite{DGLP-ldt}.

\begin{lemma}[cf.~{\cite[Lemma 5.3]{DGLP-ldt}}]\label{dt-quotient12}
Let $\Gamma$ be a connected locally $(G, s)$-distance-transitive
graph with $s\geq 2$. Let $1\neq N \lhd G$ be intransitive on
$V(\Gamma)$, and let $\mathcal{B}$ be the set of $N$-orbits on
$V(\Gamma)$. Then one of the following holds:

 \begin{itemize}
\item[(i)] $|\mathcal{B}| = 2$.

\item[(ii)] $\Gamma$ is bipartite, $\Gamma_N\cong \K_{1,r}$ with
$r\geq 2$ and $G$ is intransitive on $V(\Gamma)$.

\item[(iii)] $s=2$, $\Gamma\cong \K_{m[b]}$,   $\Gamma_N \cong \K_{m}$
with $m\geq 3$ and $b\geq 2$.

\item[(iv)] $N$ is semiregular on $V(\Gamma)$,  $\Gamma$ is a cover
of $\Gamma_N$, $|V(\Gamma_N)|<|V(\Gamma)|$ and $\Gamma_N$ is
locally $(G/N,s')$-distance-transitive where $s'=\min\{s,d(\Gamma_N)\}$.

\end{itemize}
\end{lemma}

\subsection{Cayley graphs}
Let $G$ be a finite group, and let $S$ be a subset of $G$ such that $1\notin S$ and $S=S^{-1}$. The \emph{Cayley graph} $\Cay(G,S)$ of $G$ with respect to $S$ is the graph with vertex set $G$ and edge set $\{\{g,sg\}\mid g\in G, s\in S\}$. In particular, $\Cay(G,S)$ is connected if and only if $G=\langle S\rangle$. Let $\Gamma = \Cay(G,S)$.
The group $R(G)=\{\sigma_t\mid t\in G\}$ of right translations
$\sigma_t:x\mapsto xt$ is a subgroup of $\Aut(\Gamma)$ acting regularly on
the vertex set. We may identify $G$ with $R(G)$.

\subsection{Primitive and quasiprimitive}

A transitive permutation group $G$ is said to be acting  \emph{primitively}  on a set $\Omega$  if it has only trivial blocks in $\Omega$.
If $G$ acts primitively on $\Omega$, then every nontrivial normal subgroup of $G$ is  transitive on $\Omega$.
There is a remarkable classification of finite primitive permutation
groups mainly due to M. O'Nan  and L. Scott, called the
\emph{O'Nan-Scott Theorem for primitive permutation groups}, see
\cite{LPS-1,Scott-1980}. They independently gave a classification of
finite primitive groups, and proposed their result at the ``Santa Cruz Conference on Finite Groups" in 1979. For more work on primitive groups, see~\cite{BGW-2011,LCH-circulant-2005,LS-2003}.

A graph $\Gamma$ is said to be \emph{$G$-locally-primitive} if, for each vertex $u \in V(\Gamma)$, the stabilizer $G_u$ acts primitively on the neighborhood $\Gamma(u)$.  
If $G = \Aut(\Gamma)$, we simply say that $\Gamma$ is \emph{locally-primitive}.

A transitive permutation group $G$  is \emph{quasiprimitive} if each nontrivial normal subgroup of $G$ is transitive.
This is a generalization
of primitivity as every normal subgroup of a primitive group is transitive,
but there exist quasiprimitive groups which are not primitive. For more information on quasiprimitive permutation groups, see~\cite{Praeger-4,Praeger-1}.
Praeger \cite{Praeger-4} generalized the O'Nan-Scott Theorem for primitive groups to quasiprimitive
groups and showed that a finite quasiprimitive group is one of eight distinct types: Holomorph Affine (HA), Almost Simple (AS), Twisted Wreath product (TW), Product Action (PA), Simple Diagonal (SD), Holomorph Simple (HS), Holomorph Compound (HC) and Compound Diagonal (CD).

A transitive permutation group $G$ is said to be  \emph{bi-quasiprimitive} if every nontrivial normal subgroup of $G$ has at most two orbits  and there exists one which has exactly two orbits.


\section{Quasiprimitive actions on $(G,2)$-distance-transitive graphs}

In this section, we first present a general reduction theorem for $(G,2)$-distance-transitive graphs in Subsection~\ref{Redu-Resu}, showing that the study of such graphs can be reduced to those admitting a quasiprimitive or bi-quasiprimitive group action.
Subsection~\ref{quasi-odd} focuses on $(G,2)$-distance-transitive graphs of odd order for which $G$ acts quasiprimitively on the vertex set; there we prove Theorem~\ref{2dt-qp-th1} and determine all possible vertex-quasiprimitive action types.
Building on this result, we then impose the additional assumption of $G$-local-primitivity and prove Corollary~\ref{2dt-odd-cor1}.

\subsection{Reduction result}\label{Redu-Resu}
We first establish the following reduction result for $(G,2)$-distance-transitive graphs, which follows directly from Lemma~\ref{dt-quotient12} under the additional assumption of vertex transitivity.

\begin{theo}\label{2dp-theo-1}
Let $\Gamma$ be a $(G,2)$-distance-transitive graph. Suppose that $\Gamma\ncong \K_{m[b]}$ for any $m\geq 3$ and $b\geq 2$.
Let $N$ be a normal subgroup of $G$ maximal with respect to having at least $3$ orbits.
Then the following hold:
\begin{itemize}
\item[(1)]  $N$ is semiregular on $V(\Gamma)$.
\item[(2)]  $\Gamma$ is a cover of $\Gamma_N$.
\item[(3)] Either $\Gamma_N$ is a complete $G/N$-arc-transitive graph, or $\Gamma_N$ is $(G/N,2)$-distance-transitive.
\item[(4)]  $G/N$ is quasiprimitive or bi-quasiprimitive on $V(\Gamma_N)$.
\item[(5)]  If $\Gamma$ has odd order, then $\Gamma_N$ also has odd order.
\item[(6)]  If $\Gamma$ is $G$-locally-primitive, then $\Gamma_N$ is $G/N$-locally-primitive.
\end{itemize}
\end{theo}
\proof If $N=1$, then the result is clear. Hence, in what follows, we may assume that $N\neq 1$.
Since $\Gamma$ is a $(G,2)$-distance-transitive graph, it is locally $(G,2)$-distance-transitive, and hence Lemma~\ref{dt-quotient12} applies.
Since $N$ is a normal subgroup of $G$ with at least $3$ orbits on $V(\Gamma)$, part~(i) of Lemma~\ref{dt-quotient12} does not occur. Since $G$ is transitive on $V(\Gamma)$, part~(ii) of Lemma~\ref{dt-quotient12} does not occur. Moreover, by our assumption that $\Gamma\ncong \K_{m[b]}$ for any $m\geq 3$ and $b\geq 2$, part~(iii) of Lemma~\ref{dt-quotient12} does not occur.

Therefore only part~(iv) of Lemma~\ref{dt-quotient12} can hold. Hence $N$ is semiregular on $V(\Gamma)$, $\Gamma$ is a cover of $\Gamma_N$, and $\Gamma_N$ is locally $(G/N,s')$-distance-transitive, where $s'=\min\{2,d(\Gamma_N)\}$. Since $G$ is transitive on $V(\Gamma)$, it follows that $\Gamma_N$ is $(G/N,s')$-distance-transitive.
If $d(\Gamma_N)=1$, then $\Gamma_N$ is a complete graph. Since it is $(G/N,1)$-distance-transitive, it is further $G/N$-arc-transitive. If $d(\Gamma_N)\geq 2$, then $\Gamma_N$ is $(G/N,2)$-distance-transitive. This proves~(1),~(2), and~(3).

Since $N$ is maximal among the normal subgroups of $G$ having at least $3$ orbits, every nontrivial normal subgroup of $G/N$ has at most $2$ orbits on $V(\Gamma_N)$. If some nontrivial normal subgroup of $G/N$ has exactly $2$ orbits on $V(\Gamma_N)$, then $G/N$ is bi-quasiprimitive on $V(\Gamma_N)$. On the other hand, if every nontrivial normal subgroup of $G/N$ is transitive on $V(\Gamma_N)$, then $G/N$ is quasiprimitive on $V(\Gamma_N)$.
This proves~(4).

As $N$ is semiregular on $V(\Gamma)$, each $N$-orbit on $V(\Gamma)$ has size $|N|$, and there are $|V(\Gamma_N)|$ such orbits. Therefore, $|V(\Gamma)|=|N|\cdot |V(\Gamma_N)|$.
Hence, if $\Gamma$ has odd order, then $\Gamma_N$ also has odd order. This proves~(5).

Finally, if $\Gamma$ is $G$-locally-primitive, then by \cite[Lemma~2.5]{Li-abeliancay-2008} and \cite[Lemma~4.1]{LPM-prime-2009}, $\Gamma_N$ is $G/N$-locally-primitive.
This proves~(6).
\qed

\medskip
Theorem~\ref{2dp-theo-1} provides a global analysis of the family of $(G,2)$-distance-transitive graphs. By this result, each $(G,2)$-distance-transitive graph $\Gamma\ncong \K_{m[b]}$ has a quotient graph $\Gamma_N$ corresponding to a normal subgroup $N$ of $G$ such that $\Gamma$ is a cover of $\Gamma_N$, the quotient graph $\Gamma_N$ is $(G/N,s)$-distance-transitive with $s=\min\{2,d(\Gamma_N)\}$, and $G/N$ acts quasiprimitively or bi-quasiprimitively on $V(\Gamma_N)$.
Thus, in order to classify $(G,2)$-distance-transitive graphs, we need to know all the `basic' quotient graphs, including complete graphs, for which the corresponding group acts quasiprimitively or bi-quasiprimitively on the vertex set.

\subsection{Quasiprimitive types of $(G,2)$-distance-transitive graphs of odd order}\label{quasi-odd}

This section is devoted to the study of $(G,2)$-distance-transitive graphs of odd order for which $G$ acts quasiprimitively on the vertex set, and to the proof of Theorem~\ref{2dt-qp-th1} and Corollary~\ref{2dt-odd-cor1}.

Before proving Theorem~\ref{2dt-qp-th1}, we first give examples showing that each of the three types \AS, \PA, and \HA\ can occur. These examples are constructed from Johnson graphs and Hamming graphs, so we begin by recalling their definitions.

The \emph{Johnson graph} $J(n,k)$ is defined as follows. Let $\Omega = \{1,2,\dots,n\}$, where $n \geq 3$, and let $1 \leq k \leq \left\lfloor \frac{n}{2} \right\rfloor$. The vertex set of $J(n,k)$ consists of all $k$-subsets of $\Omega$, with two vertices adjacent if and only if the corresponding $k$-subsets intersect in exactly $k-1$ elements. 
Let $\Gamma = J(n,k)$. 
By~\cite[Theorem~9.1.2]{BCN}, $\Gamma$ has diameter $k$ and valency $k(n-k)$. Moreover, $\Aut(\Gamma) \cong S_n \times \mathbb{Z}_2$ when $n=2k \geq 4$, whereas $\Aut(\Gamma) \cong S_n$ otherwise. Furthermore, $\Gamma$ is distance-transitive and, when $k\geq 2$, $2$-distance-transitive.
When $k\geq 2$, the graph $\Gamma$ has girth $3$, and hence it is not $2$-arc-transitive.

Note that $J(n,k)$ has $\binom{n}{k}=\frac{n!}{k!(n-k)!}$ vertices. For $n\geq 7$, if $n\equiv 3 \pmod 4$, then $\binom{n}{3}=\frac{n(n-1)(n-2)}{6}$ is odd. Hence, for such $n$, the Johnson graph $J(n,3)$ is a $2$-distance-transitive but not $2$-arc-transitive graph of odd order.

\begin{examp}\label{2gt-localconnect-johnson}
{\rm Let $\Gamma = J(n,k)$, where $n \geq 5$ and $2 \leq k < \frac{n}{2}$. Then its automorphism group $A := \Aut(\Gamma)$ is isomorphic to $S_n$.
Hence, for each vertex $u$, $A_u\cong S_k\times S_{n-k}$ is a maximal subgroup of $A$. Since $A$ is transitive on $V(\Gamma)$ and $n\geq 5$, it follows that $A$ is primitive of type \AS \ on $V(\Gamma)$.
}
\end{examp}

The \emph{Hamming graph} $\H(d,n)$ is the graph with vertex set $\Delta^d = \{(x_1, x_2, \ldots, x_d) \mid x_i \in \Delta\}$, where $\Delta = \{0, 1, \ldots, n - 1\}$, $d \geq 2$, and $n \geq 2$. Two vertices are adjacent if and only if they differ in exactly one coordinate.
Let $\Gamma=\H(d,n)$. Then $\Gamma$ has order $n^d$, and it is locally isomorphic to $d\K_{n-1}$.
By~\cite[Theorem~9.2.1]{BCN}, $\Gamma$ has diameter $d$ and valency $d(n-1)$. Moreover, $\Aut(\Gamma) \cong S_n \wr S_d$. Furthermore, $\Gamma$ is distance-transitive and, in particular, $2$-distance-transitive.
If $n\geq 3$, then $\Gamma$ has girth $3$, and hence it is not $2$-arc-transitive.

\begin{examp}\label{2gt-qp-examphamm1}
{\rm  
Let $\Gamma=\H(d,n)$, where $d\geq 2$ and $n$ is an odd prime power. Then $|V(\Gamma)|=n^d$ is odd. Set $A:=\Aut(\Gamma)\cong S_n \wr S_d$.
If $n=3$, then $A$ acts primitively of type \HA \ on $V(\Gamma)$.
If $n \geq 5$, then $A$ acts primitively of type \PA \ on $V(\Gamma)$.
}
\end{examp}

We now prove Theorem~\ref{2dt-qp-th1}, which determines all possible quasiprimitive action types.

\medskip
\noindent {\bf Proof of Theorem~\ref{2dt-qp-th1}.}
Let $\Gamma$ be a $(G,2)$-distance-transitive graph of odd order, and suppose that $G$ is quasiprimitive on $V(\Gamma)$ of type $X$.
Then, by \cite{Praeger-4}, $X$ is one of the following eight types: HA, AS, HS, HC, SD, CD, PA, and TW.

Suppose that $X$ is of type \TW, \HS, or \HC. In each of these cases, $G$ contains a minimal normal subgroup $M$ acting regularly on $V(\Gamma)$, where $M\cong T^m$ for some integer $m\geq 1$ and some nonabelian simple group $T$. Thus $|V(\Gamma)|=|M|=|T|^m$. However, by the Feit--Thompson theorem, every nonabelian simple group has even order. It follows that $|V(\Gamma)|$ is even, contradicting the assumption that $|V(\Gamma)|$ is odd.

If  $X$ is of type \SD\, then $\soc(G)\cong T^k$ with $k\geq 2$, where $T$ is a nonabelian simple group. In this case, $T^{k-1}$ acts regularly on $V(\Gamma)$. 
Therefore, $|V(\Gamma)|=|T^{k-1}|=|T|^{k-1}$.
However, every nonabelian simple group has even order, and hence $|V(\Gamma)|$ is even, contrary to the assumption that $|V(\Gamma)|$ is odd.

Assume that $X$ is of type \CD.
Then $V(\Gamma)=\Delta^l$ and $\soc(G)=T^k\leq G\leq N\wr S_l\leq \Sym(\Delta)\wr S_l$, for some divisor $l$ of $k$ and some nonabelian simple group $T$, where $l\geq 2$, $k/l\geq 2$, $N\leq \Sym(\Delta)$, $\soc(N)=T^{k/l}$, and $N$ is quasiprimitive of SD type on $\Delta$. By the previous paragraph, $|\Delta|$ is even, and hence $|V(\Gamma)|=|\Delta|^l$ is even, a contradiction.

Therefore, $X$ is one of the following three types: \AS, \PA, and \HA.
Moreover, by Examples~\ref{2gt-localconnect-johnson} and~\ref{2gt-qp-examphamm1}, there are infinitely many examples for each type. 

This completes the proof.
\qed

Finally, we further assume that $\Gamma$ is $G$-locally-primitive and prove Corollary~\ref{2dt-odd-cor1}.

\medskip
\noindent {\bf Proof of Corollary~\ref{2dt-odd-cor1}.}
Assume that $G$ acts quasiprimitively on $V(\Gamma)$ of type $X$. Then it follows from Theorem~\ref{2dt-qp-th1} that $X$ is one of the following three types: \AS, \PA, and \HA. Therefore, to prove Corollary~\ref{2dt-odd-cor1}, it suffices to show that $X$ cannot be of type \HA.

Suppose that $X$ is of type \HA.
Then $G$ has a regular normal subgroup $N\cong \mathbb{Z}_p^d$ for some prime $p$ and integer $d$.
Let $u\in V(\Gamma)$.
We can identify $V(\Gamma)$ with $N$ such that $u$ corresponds to the identity element of $N$; that is, $\Gamma\cong \Cay(N,S)$, where
$S = S^{-1}=\Gamma(u)\subset N\setminus \{u\}$ and $N=\langle S\rangle$.

Let $v\in \Gamma(u)$. Then $\{v,v^{-1}\}\subseteq \Gamma(u)$.
Since $N$ is normal and regular in $G$, the conjugation action embeds $G_u$ into $\Aut(N)$. Hence $\{v,v^{-1}\}$ is a block of the $G_u$-action on $\Gamma(u)$.
As $G_u$ is primitive on $\Gamma(u)$, the set
$\{v,v^{-1}\}$ must be a trivial block. Since the valency of $\Gamma$ is at least $3$, we have $|\{v,v^{-1}\}|=1$, and so $v=v^{-1}$. Therefore, $p=2$, and hence $|V(\Gamma)|=2^d$, a contradiction.
Therefore $X$ is not of type \HA.

This proves Corollary~\ref{2dt-odd-cor1}.
\qed


\section{$2$-Distance-transitive graphs of valency $p+1$ }
In this section, we study $2$-distance-transitive graphs of valency $p+1$, where $p$ is a prime, and prove Theorem~\ref{2dt-theo-prime1}.

Let $\Gamma$ be a connected $G$-arc-transitive graph, and let $(u,v)$ be an arc. Recall that $\Gamma$ is $(G,2)$-geodesic-transitive if and only if $G_{uv}$ is transitive on $\Gamma_2(u)\cap \Gamma(v)$.

\begin{lemma}\label{2dt-2gt-1}
Let $\Gamma$ be a connected $(G,2)$-distance-transitive graph of valency
$k\geq 2$. If $(b_1,k)=1$ and $c_2\mid k$, then $\Gamma$ is
$(G,2)$-geodesic-transitive. In particular, if $(b_1,k)=1$ and
$(b_1,c_2)=1$, then $\Gamma$ is $(G,2)$-geodesic-transitive.
\end{lemma}
\proof Let $(u,v,w)$ be a $2$-geodesic of $\Gamma$. Since $\Gamma$ is $(G,2)$-distance-transitive, $G_u$ is transitive on both $\Gamma(u)$ and $\Gamma_2(u)$. Counting the edges between $\Gamma(u)$ and $\Gamma_2(u)$ gives $kb_1=c_2|\Gamma_2(u)|$. Since $G_u$ is transitive on $\Gamma_2(u)$, we have $|\Gamma_2(u)|=|G_u:G_{uw}|$. As $c_2\mid k$, it follows that $b_1$ divides $\frac{b_1k}{c_2}=|\Gamma_2(u)|=|G_u:G_{uw}|.$

Moreover,
$$
|G_u:G_{uw}|\cdot |G_{uw}:G_{uwv}|
=|G_u:G_{uv}|\cdot |G_{uv}:G_{uvw}|
=k|G_{uv}:G_{uvw}|.
$$
Thus $b_1$ divides $k|G_{uv}:G_{uvw}|$. Since $(k,b_1)=1$, $b_1$ divides $|G_{uv}:G_{uvw}|$. On the other hand, $|G_{uv}:G_{uvw}|\leq |\Gamma_2(u)\cap \Gamma(v)|=b_1.$
Hence $|G_{uv}:G_{uvw}|=b_1$, and so $G_{uv}$ is transitive on $\Gamma_2(u)\cap \Gamma(v)$. Therefore, $\Gamma$ is $(G,2)$-geodesic-transitive.

Finally, if $(b_1,c_2)=1$, then the identity $kb_1=c_2|\Gamma_2(u)|$ implies that $c_2\mid k$, and hence the first part applies.
\qed

\begin{lemma}\label{2dt-girth4-1}
Let $\Gamma$ be a connected $(G,2)$-distance-transitive
graph of valency $k\geq 2$. Suppose that $(c_2,k-1)=1$.
Then either $\Gamma$ has girth $3$ or $\Gamma$ is $(G,2)$-arc-transitive.
\end{lemma}
\proof Suppose that $\Gamma$ has girth at least $4$. Then $b_1=k-1$.
Since $(c_2,k-1)=1$, we have $(c_2,b_1)=1$. Moreover,
$(b_1,k)=(k-1,k)=1$. As $\Gamma$ is $(G,2)$-distance-transitive,
Lemma~\ref{2dt-2gt-1} implies that $\Gamma$ is
$(G,2)$-geodesic-transitive. Since $\Gamma$ has girth at least $4$,
every $2$-arc in $\Gamma$ is a $2$-geodesic. Hence, $\Gamma$ is
$(G,2)$-arc-transitive.
\qed

\medskip

\noindent {\bf Proof of Theorem~\ref{2dt-theo-prime1}.}
The necessity is obvious. We now prove the sufficiency. 
Assume that $\Gamma$ has girth at least $4$. Then $b_1=p$ and $1\leq c_2\leq p+1$. If $c_2\neq p$, then $(c_2,p)=1$, so Lemma~\ref{2dt-girth4-1} implies that $\Gamma$ is $2$-arc-transitive. If $c_2=p$, then $\Gamma$ has girth $4$, and by~\cite[Lemma~4.4(ii)]{CJS-2dt}, $\Gamma\cong\K_{p+2,p+2}-(p+2)\K_2$, which is $2$-arc-transitive.
\qed


\section{$2$-Distance-transitive graphs of small valency}

In this section, we study locally-primitive $2$-distance-transitive graphs of valency at most $8$. Our main objective is to prove Theorem~\ref{2dt-localp-th1}.
The following fact will be used repeatedly in this section.

\begin{lemma}\label{knn} 
The following statements hold:
\begin{itemize}
\item[(1)] The complement of the $(2 \times n)$-grid is isomorphic to $\K_{n,n} - n\K_2$, and is $2$-arc-transitive.
\item[(2)] The complete bipartite graph $\K_{n,n}$ is $2$-arc-transitive.
\end{itemize}
\end{lemma}

\subsection{Restrictions on parameters}

First, we present several lemmas that impose strong restrictions on the graphs. 

\begin{lemma}\label{2dt-parameter} 
Let $\Gamma$ be a connected $2$-distance-transitive graph of valency $k \geq 3$. Then the following statements hold:
\begin{itemize}
    \item[(1)] If $k$ is odd, then both $a_1$ and $b_1$ are even.
    \item[(2)] If $b_1 = 1$, then $k$ is even and $\Gamma \cong \K_{\frac{k+2}{2}[2]}$.
    \item[(3)] If $\Gamma$ has diameter $d \geq 3$ and $b_1 = c_2$, then it is a Taylor graph. Moreover, every $2$-distance-transitive Taylor graph is distance-transitive.
    \item[(4)] If $\Gamma$ has diameter $d \geq 3$ and $b_1 = c_2 = k-1$, then $\Gamma \cong \K_{k+1,k+1}-(k+1)\K_2$.
    \item[(5)] If $\Gamma$ has parameters $b_1 = k-1$ and $c_2 = k$, then $\Gamma \cong \K_{k,k}$.

\end{itemize}
\end{lemma}

\proof
Item (1)  follows from applying the handshaking lemma to the induced subgraph $[\Gamma(u)]$ for a vertex $u \in V(\Gamma)$.

Assume that $b_1 = 1$. Then, by (1), $k$ is even, and hence $k \geq 4$. Let $(u, v, w)$ be a $2$-geodesic in $\Gamma$. Since $b_1 = 1$, we have $a_1 = k - 2$, and thus $ |\Gamma(u) \cap \Gamma(v)| = |\Gamma(v) \cap \Gamma(w)| = k - 2$.
Therefore, $\Gamma(v) = \{u, w\} \cup (\Gamma(u) \cap \Gamma(v))$, and hence
$\Gamma(u)\cap\Gamma(v)=\Gamma(v)\cap\Gamma(w)$, which implies
$c_2 = |\Gamma(u)\cap\Gamma(w)| \geq k-2$.
By counting the edges between $\Gamma(u)$ and $\Gamma_2(u)$, we obtain $k b_1 = |\Gamma_2(u)| c_2$, and hence $k \geq |\Gamma_2(u)|(k-2)$. 
Assume that $|\Gamma_2(u)| \geq 2$, then $k=4$, $c_2=2$, and $|\Gamma_2(u)|=2$.
If $d \geq 3$, then by \cite[Proposition~4.1.6]{BCN} we have $c_2 \leq b_1$, contradicting $c_2 = 2$ and $b_1 = 1$. If $d = 2$, then $a_2 = k - c_2 = 2$, which contradicts the fact that $|\Gamma_2(u)| = 2$. It follows that $|\Gamma_2(u)| = 1$.
Hence, $\Gamma_2(u) = \{w\}$, and consequently $|\Gamma(u) \cap \Gamma(w)| = k$.  
It is straightforward to verify that $|V(\Gamma)| = k + 2$ and $\Gamma \cong \K_{\frac{k+2}{2}[2]}$. This proves (2).

We now proceed to prove (3).  
Assume that $b_1 = c_2$. Then we have $k=a_1+b_1+c_1=a_1+c_2+1$, and by Theorem~\ref{Thm1.5.5}, $\Gamma$ is a Taylor graph.
Let $u \in V(\Gamma)$. Then $V(\Gamma)=\{u\}\cup \Gamma_{1}(u)\cup \Gamma_{2}(u)\cup \Gamma_{3}(u)$, where $|\Gamma_{2}(u)|=k$ and $|\Gamma_{3}(u)|=1$.
Therefore, for each vertex $u \in V(\Gamma)$, there exists a unique vertex $v \in V(\Gamma)$ such that $d_{\Gamma}(u, v) = 3$.  
Since $\Gamma$ is $2$-distance-transitive, it must also be $3$-distance-transitive, and hence distance-transitive. This proves (3).

Next, assume that $\Gamma$ has diameter $d \geq 3$ and satisfies $b_1 = c_2 = k - 1$. Then, by (3), $\Gamma$ is a Taylor graph with intersection array $(k, k - 1, 1; 1, k - 1, k)$.  
Assertion~(4) follows directly from \cite[Corollary~1.5.4]{BCN}.

Finally, to prove (5), assume that $c_2 = k$. Then $a_2 = b_2 = 0$, so $\Gamma$ has diameter $2$.  
Let $u \in V(\Gamma)$. Then $V(\Gamma) = \{u\} \cup \Gamma(u) \cup \Gamma_2(u)$.  
By counting the edges between $\Gamma(u)$ and $\Gamma_2(u)$, we find that $|\Gamma_2(u)| = k - 1$.  
It follows that $\Gamma \cong \K_{k,k}$. This proves (5).
\qed

\begin{lemma}\label{local-p}
Let $\Gamma$ be a connected $2$-distance-transitive graph of valency $k \geq 3$. Then the following statements hold:

\begin{itemize}
    \item[(1)] If $a_1 \geq 1$ and the local graph $[\Gamma(u)]$ is disconnected for some vertex $u \in V(\Gamma)$, then $\Gamma$ is not locally-primitive.
    \item[(2)] If $a_1 = 1$, then $\Gamma$ is not locally-primitive.
    \item[(3)] If $a_1 \geq c_2 = 1$, then $\Gamma$ is not locally-primitive.
    \item[(4)] If $c_2 = 1$, then either $\Gamma$ is $2$-arc-transitive, or it is neither $2$-arc-transitive nor locally-primitive.
\end{itemize}
\end{lemma}

\proof
Let $u \in V(\Gamma)$. Since $[\Gamma(u)]$ is an $a_1$-regular graph on $k \geq 3$ vertices, if $a_1 \geq 1$ then each connected component of $[\Gamma(u)]$ contains at least two vertices. Therefore, (1) holds.

If $a_1 = 1$, then $[\Gamma(u)]$ consists of disjoint edges, so (2) follows from (1).

To prove (3), let $(u, v, w)$ be a $2$-geodesic in $\Gamma$. Since $c_2 = 1$, we have $\Gamma(u) \cap \Gamma(w) = \{v\}$, and hence $w$ is not adjacent to any vertex in $\Gamma(u) \cap \Gamma(v)$. As $w$ can be any vertex in $\Gamma_2(u) \cap \Gamma(v)$, it follows that $[\Gamma(v)]$ is disconnected, with each component containing at least $a_1 + 1 \geq 2$ vertices. Hence, $\Gamma$ is not locally-primitive.

Assume that $c_2 = 1$. If $a_1 \geq 1$, then by (3), $\Gamma$ is not locally-primitive. If $a_1 = 0$, then $\Gamma$ has girth at least $5$, and hence is $2$-arc-transitive. Therefore, (4) holds.
\qed

\begin{lemma}\label{primitive-small-local}
Let $\Delta$ be a graph on $4$, $6$, or $8$ vertices. If a subgroup of $\Aut(\Delta)$ acts primitively on $V(\Delta)$, then $\Delta$ is complete or empty.
\end{lemma}
\proof
Let $A\leq\Aut(\Delta)$ act primitively on $V(\Delta)$, and suppose that $\Delta$ is neither complete nor empty. Since $A$ is transitive, $\Delta$ is regular. Replacing $\Delta$ by its complement if necessary, we may assume that its valency $r$ satisfies $1\leq r\leq (|V(\Delta)|-1)/2$. The connected components form an $A$-invariant partition, so $\Delta$ is connected. Moreover, $\Delta$ cannot be bipartite, since its two bipartite halves would form a nontrivial $A$-invariant partition. Thus $r\neq 1$, and $r\neq 2$ because a connected $2$-regular graph of even order is an even cycle. It follows that $|V(\Delta)|=8$ and $r=3$.

Let $t$ be the number of triangles containing a vertex of $\Delta$; this is independent of the vertex. Counting vertex--triangle incidences gives $3\mid 8t$. Since $0\leq t\leq 3$, we have $t=0$ or $3$. If $t=3$, then a vertex together with its three neighbors forms a component isomorphic to $\K_4$, contrary to connectedness. Hence $\Delta$ is triangle-free.

Fix $x\in V(\Delta)$, write $\Delta(x)=\{a,b,c\}$, and let $Y=V(\Delta)\setminus(\{x\}\cup\Delta(x))$. There are six edges between $\Delta(x)$ and $Y$, so $[Y]$ has three edges. Since $[Y]$ is triangle-free on four vertices, it is either a star $\K_{1,3}$ or a path on four vertices. In the star case, $\{x\}$ together with the three leaves is one bipartite half of $\Delta$, a contradiction.

In the path case, write $Y=\{y_1,y_2,y_3,y_4\}$ with path $y_1y_2y_3y_4$. Each endpoint has two neighbors in $\{a,b,c\}$, and each internal vertex has one. The endpoints cannot have the same two neighbors: otherwise both internal vertices would be adjacent to the third vertex of $\{a,b,c\}$, creating a triangle. Relabelling $a,b,c$ if necessary, we therefore have $\Delta(y_1)\cap\Delta(x)=\{a,b\}$ and $\Delta(y_4)\cap\Delta(x)=\{a,c\}$. Triangle-freeness then forces $y_2\sim c$ and $y_3\sim b$. These adjacencies determine $\Delta$: it consists of the cycle $x,b,y_3,y_4,a,y_1,y_2,c,x$ and the four edges joining opposite vertices on this cycle. Each of these four edges lies in two $4$-cycles, whereas every cycle edge lies in exactly one. Hence the four opposite pairs form an $\Aut(\Delta)$-invariant partition of $V(\Delta)$, contradicting the primitivity of $A$.
\qed

\subsection{Discussion on valency $7$}

We begin the proof of Theorem~\ref{2dt-localp-th1} by treating the case $k = 7$.

For graphs of valency $7$, we first consider the following special case.

\begin{lemma}\label{741127-not-exist}
There is no $2$-distance-transitive graph of valency $7$ with $a_1=c_2=2$.
\end{lemma}
\proof
Suppose, to the contrary, that such a graph $\Gamma$ exists. Since $a_1=c_2=2$, the graph $\Gamma$ is a $(0,2)$-graph; moreover, it is non-bipartite because $a_1\neq 0$. Since $b_1=7-a_1-1=4$, counting the edges between $\Gamma(v)$ and $\Gamma_2(v)$ gives $k_2=7b_1/c_2=14$.

By the classification in~\cite[Table~2]{Brouwer-2006}, the only vertex-transitive $(0,2)$-graph satisfying these conditions is the graph listed in line~12. For each vertex $v$, the distance partition has cell sizes $k_0=1$, $k_1=7$, $k_2=14$, and $k_3=2$. This graph is distance-regular with intersection array $\{7,4,1;1,2,7\}$.

By~\cite[p.~386]{BCN}, the automorphism group $G=\Aut(\Gamma)$ acts transitively on $V(\Gamma)$ with rank $6$. Thus, the point stabilizer $G_v$ has exactly six orbits on $V(\Gamma)$. Since $\Gamma$ is $2$-distance-transitive, both $\Gamma_1(v)$ and $\Gamma_2(v)$ are single $G_v$-orbits. Together with $\{v\}$, they account for three of the six orbits, so the remaining three orbits must be contained in $\Gamma_3(v)$. This is impossible because $|\Gamma_3(v)|=2$.
\qed

\begin{lemma}\label{2dt-localp-val7}
Let $\Gamma$ be a $2$-distance-transitive graph of valency $7$.
If $\Gamma$ is locally-primitive, then $\Gamma$ is $2$-arc-transitive.
\end{lemma}
\proof
We consider the possible values of $a_1$.

\noindent
\textbf{Case 1: $a_1\geq 1$.}
The local graph $[\Gamma(u)]$ is an $a_1$-regular graph of order $7$. By Lemma~\ref{2dt-parameter}(1), both $a_1$ and $b_1$ are even. Therefore, $(a_1,b_1)\in\{(2,4),(4,2)\}$.

Assume first that $d=2$. Then $\Gamma$ is a strongly regular graph. If $(a_1,b_1)=(2,4)$, then $28=kb_1=c_2k_2$. Since $c_2\leq 7$, we have $k_2\in\{4,7,14,28\}$, and hence $|V(\Gamma)|\in\{12,15,22,36\}$. If $(a_1,b_1)=(4,2)$, then $14=kb_1=c_2k_2$, so $k_2\in\{2,7,14\}$ and $|V(\Gamma)|\in\{10,15,22\}$. However, by~\cite[Chapter~12]{BM-2022}, there is no strongly regular graph of valency $7$ whose order belongs to $\{10,12,15,22,36\}$.

Now assume that $d\geq 3$. By Lemma~\ref{b1geqk}, we have $b_1\geq \frac{1}{3}(k+1)=\frac{8}{3}$, and hence $(a_1,b_1)=(2,4)$. By Theorem~\ref{Thm1.5.5}, we have $c_2\leq k-a_1-1=4$. Moreover, $28=kb_1=c_2k_2$, so $c_2\mid 28$ and therefore $c_2\in\{1,2,4\}$. If $c_2=1$, then Lemma~\ref{local-p}(3) implies that $\Gamma$ is not locally-primitive, a contradiction. If $c_2=2$, then Lemma~\ref{741127-not-exist} yields a contradiction. Finally, if $c_2=4$, then Lemma~\ref{2dt-parameter}(3) implies that $\Gamma$ is a distance-transitive Taylor graph with intersection array $\{7,4,1;1,4,7\}$. However, by~\cite[Theorem~7.5.3]{BCN}, no such distance-transitive graph exists, again a contradiction.

\noindent
\textbf{Case 2: $a_1=0$.}
Then $b_1=6$.
Let $(u,v,w)$ be a $2$-geodesic, and let $k_2 = |\Gamma_2(u)|$. By counting the edges between $\Gamma_1(u)$ and $\Gamma_2(u)$, we obtain $42 = k b_1 = k_2 c_2$.

If $c_2 = 7$, then by Lemma~\ref{2dt-parameter}(5), we have $\Gamma \cong \K_{7,7}$, which is $2$-arc-transitive.

Assume that $c_2 = 6$, then $k_2 = 7$. If the diameter $d = 2$, then $[\Gamma_2(u)]$ is a regular graph of valency $1$ and order $7$, which is impossible. Consequently, $d \geq 3$. By Lemma~\ref{2dt-parameter}(4), we have $\Gamma \cong \K_{8,8} - 8\K_2$, which is $2$-arc-transitive.

Suppose that $c_2 \leq 5$. Since $42 = k_2 c_2$, it follows that $7 \mid k_2$ and $c_2 \mid 6$, and hence $c_2 \leq 3$. 
If $\Gamma$ has girth at least $5$, then it is $2$-arc-transitive. Therefore, we may assume that $\Gamma$ has girth $4$, and hence $c_2\geq 2$. It follows that $c_2=2$ or $3$.

Assume that $c_2 = 3$. Then, by \cite[Corollary 4.8]{JKL-2025-1}, either $d=3$ and $|V(\Gamma)|=30$, or $d=4$ and $|V(\Gamma)|=32$. According to \cite{Conder-website}, for each of these two orders there is a unique connected arc-transitive graph of valency $7$, and both graphs are $2$-arc-transitive. Thus $\Gamma$ is $2$-arc-transitive.

Finally, assume that $c_2 = 2$. Let $G := \Aut(\Gamma)$. Then $\Gamma$ is an amply regular graph with $a_1 = 0$, $c_2 = 2$, and $|\Gamma_2(u)| = 21$. 
In particular, $\Gamma$ is a $(0,2)$-graph.
Since $\Gamma$ is $2$-distance-transitive, it follows that both $|\Gamma(u)|$ and $|\Gamma_2(u)|$ divide $|G_u|$.

Suppose that $\Gamma$ is bipartite. Consider the graphs of valency $7$ listed in \cite[Table~1]{Brouwer-2006} whose automorphism group is transitive and for which $|\Gamma_2(u)| = 21$ divides $|G|$. Then $\Gamma$ must be one of the graphs listed in lines~10, 29, 39, or~40 of \cite[Table~1]{Brouwer-2006}.
If $\Gamma$ corresponds to the graph listed in line~39, then $|G_u| = |G|/|V(\Gamma)| = 16128/112 = 144$. However, $|\Gamma_2(u)| = 21$ does not divide $|G_u|$, and hence this case does not occur.
If $\Gamma$ corresponds to the graph listed in line~10, 29, or~40, then $|G_u| = 42$, $168$, or $5040$, respectively. In each case, by \cite{LLW-2016} or \cite[Lemma~2.5]{PHW-2021}, $\Gamma$ is $2$-arc-transitive.

Suppose that $\Gamma$ is not bipartite. Consider the graphs of valency $7$ listed in \cite[Table~2]{Brouwer-2006} whose automorphism group is transitive and for which $|\Gamma_2(u)| = 21$ divides $|G|$. Then $\Gamma$ must be the graph listed in line~51 of \cite[Table~2]{Brouwer-2006}. Hence $|V(\Gamma)| = 64$ and $|G| = 322560$, so $|G_u| = 5040$. By \cite{LLW-2016} or \cite[Lemma~2.5]{PHW-2021}, $\Gamma$ is $2$-arc-transitive.

This completes the proof.
\qed

\subsection{Proof of Theorem~\ref{2dt-localp-th1}}

\noindent {\bf Proof of Theorem~\ref{2dt-localp-th1}.}
Let $\Gamma$ be a $2$-distance-transitive graph of valency $k$, where $2 \leq k \leq 8$. Suppose that $\Gamma$ is locally-primitive. Clearly, $\Gamma$ is not complete, and hence has diameter $d \geq 2$. If $k\in\{4,6,8\}$, then Lemma~\ref{primitive-small-local} applied to the local graph $[\Gamma(u)]$ shows that it is complete or empty. Since $\Gamma$ is connected and vertex-transitive with diameter at least $2$, its local graph cannot be complete. Thus $a_1=0$, and $\Gamma$ has girth at least $4$. As $k-1\in\{3,5,7\}$ is prime, Theorem~\ref{2dt-theo-prime1} implies that $\Gamma$ is $2$-arc-transitive. By Lemma~\ref{2dt-localp-val7}, we may therefore assume that $k\in\{2,3,5\}$.

If $k = 2$, then $\Gamma$ is the cycle graph $C_n$ with $n \geq 4$, and so it is $2$-arc-transitive.

Let $k = 3$. If $\Gamma$ is not $2$-arc-transitive, then by \cite[Theorem 1.3]{CJS-2dt}, the only candidate is the complement of the $(2 \times 4)$-grid. This graph is $2$-arc-transitive, as shown by Lemma~\ref{knn}.

Let $k = 5$. If $\Gamma$ is not $2$-arc-transitive, then according to \cite[Theorem 1.3]{CJS-2dt}, the only possibilities are the icosahedron or the complement of the $(2 \times 6)$-grid. Note that the complement of the $(2 \times 6)$-grid is $2$-arc-transitive. Thus $\Gamma$ must be the icosahedron.

This completes the proof of Theorem~\ref{2dt-localp-th1}.
\qed

\section{Restrictions on square-free orders in locally-primitive graphs}
This section is devoted to the proof of Theorem~\ref{2dt-lp-sf}. Our main strategy is to analyze the Fitting subgroup of the soluble group $G$. 

The characteristic subgroup  of a group $G$  generated by all the nilpotent normal subgroups of $G$ is called the \emph{Fitting subgroup} of $G$.

The following are  some technical results which will be used to prove Theorem~\ref{2dt-lp-sf}.

\begin{lemma}\label{2gt-sqf-soluble-1}
Let $\Gamma$ be a $(G,2)$-distance-transitive graph of square-free order, and assume that $G$ is soluble.
Let $F$ be the Fitting subgroup of $G$.
Then one of the following holds:
\begin{itemize}
\item[(1)] $|F|$ is square-free.
\item[(2)] $\Gamma \cong \K_{n[b]}$, where $n \geq 3$, $b \geq 2$, and $nb$ is square-free.
\item[(3)] $\Gamma \cong \K_{p,p}$, where $p$ is an odd prime and $|F| = p^i m$ with $i \geq 2$ and $m$ square-free.
\end{itemize}
\end{lemma}
\proof Since $G$ is soluble, $F$ is nontrivial. By definition, $F$ is a nilpotent normal subgroup of $G$. In a finite nilpotent group, every Sylow $p$-subgroup is normal and hence unique.
Let $P$ be a Sylow $p$-subgroup of $F$. Since $P$ is the unique Sylow $p$-subgroup of $F$, it is invariant under all automorphisms of $F$, and hence $P$ is a characteristic subgroup of $F$.
Since a characteristic subgroup of a normal subgroup is normal in the whole group, we have $P \lhd G$.
Thus every Sylow subgroup of $F$ is normal in $G$.

If $|F|$ is square-free, then the first conclusion of the lemma holds. Thus, we may assume that $|F|$ is not square-free. 
Then $F$ has a Sylow $p$-subgroup $P$ such that  $|P|= p^i$ where $i\geq 2$.
By the above discussion, we have $P \lhd G$.

Suppose that $P$ has at least three orbits on $V(\Gamma)$.
Since $\Gamma$ is $(G,2)$-distance-transitive, $G$ is transitive on $V(\Gamma)$.
Applying Lemma~\ref{dt-quotient12} with the normal subgroup $P$, 
our assumption on the number of orbits rules out Case~(i), and the transitivity of $G$ rules out Case~(ii).
Consequently, either Case~(iii) or Case~(iv) of Lemma~\ref{dt-quotient12} occurs.
If Case~(iii) occurs, then $\Gamma \cong \K_{n[b]}$, where $n \geq 3$, $b \geq 2$, and $nb=|V(\Gamma)|$ is square-free, which gives conclusion~(2).
If Case~(iv) occurs, then $P$ acts semiregularly on $V(\Gamma)$.
This implies that every orbit of $P$ has length exactly $|P|=p^i$.
Consequently, $p^i$ divides $|V(\Gamma)|$, contradicting the assumption that $|V(\Gamma)|$ is square-free.

If $P$ is transitive on $V(\Gamma)$, then $|V(\Gamma)|$ divides $|P|$. Since $|V(\Gamma)|$ is square-free,
we have $|V(\Gamma)|=p$. Thus $P\leq \Sym(p)$, and so $|P|=p^i$ with $i\geq 2$ divides $p!$, which is impossible.

Finally, assume that $P$ has exactly 2 orbits on $V(\Gamma)$, say $U_0$ and $U_1$. 
Since $P$ is normal in $G$, the orbit set $\{U_0, U_1\}$ forms a $G$-invariant block system. 
As $\Gamma$ is connected and $G$-arc-transitive, 
the induced subgraph $[U_j]$ contains no edge of $\Gamma$ for $j=0,1$.
Hence, $\Gamma$ is a bipartite graph with biparts $U_0$ and $U_1$.
Since $\Gamma$ is vertex-transitive, it is a regular graph, which implies $|U_0| = |U_1|$. 
As $P$ acts transitively on $U_j$, it follows that $|U_j|$ divides $|P|$. Hence $|U_j|$ is a power of $p$ for $j=0,1$.
Since $|U_0| = \frac{|V(\Gamma)|}{2}$ is square-free, we have $|U_0| = |U_1| = p$.

If $P$ acts faithfully on some $U_j$, then $P\leq \Sym(p)$, and hence $p^i=|P|$ divides $p!$, where $i\geq 2$, a contradiction.
Thus $P$ is not faithful on $U_0$ and $U_1$.
Let $P_{(U_0)}$ be the kernel of $P$ on $U_0$.
Then $|P_{(U_0)}|=p^r$ for some $1\leq r<i$. Let $g\in P_{(U_0)}$ be an element of order $p$. Then $g$
fixes each vertex of $U_0$ pointwise. Since $g$ has order $p$ and $|U_1|=p$,
the cyclic group $\langle g \rangle$ either acts transitively on $U_1$ or acts trivially on $U_1$.
If the latter case holds, then $g$ fixes each vertex of $V(\Gamma)$ pointwise, and so $g=1$, a contradiction. Thus
$g$ is transitive on $U_1$, which forces every vertex of $U_0$ to be adjacent to every vertex of $U_1$. Therefore $\Gamma\cong \K_{p,p}$.
For each prime $q$ such that $q^i$ divides $|F|$ for some $i\geq 2$, the above discussion implies that $\Gamma\cong \K_{q,q}$. Therefore $q=p$, and hence $|F|/p^i$ is square-free.
Thus $|F|=p^i m$, where $i\geq 2$ and $m$ is square-free. This gives conclusion~(3).
\qed

\begin{lemma}\label{2dt-sqf-quot-1}
Let $\Gamma$ be a $(G,2)$-distance-transitive graph of square-free order, and assume that $G$ is soluble.
Let $F$ be the Fitting subgroup of $G$. If $|F|$ is square-free, then $F = C_G(F)$ is cyclic, and one of the following holds:
\begin{itemize}
    \item[(1)] $\Gamma \cong \K_{m[b]}$, where $m \geq 3$, $b \geq 2$, and $mb$ is square-free.  
    \item[(2)] $F$ acts regularly on $V(\Gamma)$.
\end{itemize}
\end{lemma}
\proof
Since $F$ is nilpotent and has square-free order, each Sylow subgroup of $F$ is normal in $F$ and has prime order, and the orders of distinct Sylow subgroups are coprime. Hence $F$ is cyclic. In particular, $\Aut(F)$ is abelian, and $F\leq C_G(F)$.
On the other hand, since $G$ is soluble and $F$ is the Fitting subgroup of $G$, we have $C_G(F)\leq F$.
Therefore $C_G(F)=F$, and hence $G/F=N_G(F)/C_G(F)\leq \Aut(F)$ is abelian.
Since $F$ is a cyclic normal subgroup of $G$ and $G$ is transitive on $V(\Gamma)$, it follows that $F$ acts semiregularly on $V(\Gamma)$.

\noindent
\textbf{Case 1: $F$ has at least three orbits on $V(\Gamma)$.}

Since $\Gamma$ is $(G,2)$-distance-transitive, $F\lhd G$, and $F$ has at least three orbits on $V(\Gamma)$, Lemma~\ref{dt-quotient12} applies. In particular, only Case~(iii) or Case~(iv) of Lemma~\ref{dt-quotient12} can occur.
If Case~(iii) occurs, then $\Gamma\cong \K_{m[b]}$ with $m\geq 3$ and $b\geq 2$, giving the first alternative.

Suppose that Case~(iv) occurs. Since $G$ is transitive on $V(\Gamma)$, the group $G/F$ is transitive on $V(\Gamma_F)$. Moreover, $\Gamma$ is a cover of $\Gamma_F$, and $\Gamma_F$ is $(G/F,s)$-distance-transitive, where $s=\min\{2,d(\Gamma_F)\}$.
In particular, $G/F$ is transitive on the arcs of $\Gamma_F$.
Since $G/F$ is abelian, the faithful permutation group it induces on $V(\Gamma_F)$ is regular.
As $\Gamma_F$ has at least three vertices and $\Gamma$ is a cover of $\Gamma_F$, its valency is at least $2$. Thus $G/F$ cannot be transitive on the arcs of $\Gamma_F$, a contradiction.

\noindent
\textbf{Case 2: $F$ has at most two orbits on $V(\Gamma)$.}

If $F$ is transitive on $V(\Gamma)$, then, as $F$ acts semiregularly on $V(\Gamma)$, the second alternative holds.
Therefore, to prove Lemma~\ref{2dt-sqf-quot-1}, it remains to show that $F$ cannot have two orbits on $V(\Gamma)$.
For the remainder of the proof, assume that $F$ has two orbits on $V(\Gamma)$, and we shall derive a contradiction.

Since $\Gamma$ is $G$-arc-transitive and $F$ is a normal subgroup of $G$, each $F$-orbit contains no edge of $\Gamma$. Hence $\Gamma$ is bipartite, and the two $F$-orbits are the two bipartite halves of $\Gamma$, say $U_1$ and $U_2$.
Let $m=|F|\geq 2$, and let $k\geq 2$ be the valency of $\Gamma$. Since $F$ is semiregular, it is regular on each $U_j$, so $|V(\Gamma)|=2m$. For every vertex $u$, $G_u\cap F=1$, and hence $G_u$ embeds in the abelian group $G/F$.

Let $(u,v)$ be an arc. Since $G_u$ is abelian and transitive on $\Gamma(u)$, every element of $G_{uv}$ fixes each vertex $v^g$ with $g\in G_u$, and thus fixes $\Gamma(u)$ pointwise. Similarly, it fixes $\Gamma(v)$ pointwise. Whenever $G_{uv}$ fixes adjacent vertices $x,y$, it lies in $G_{xy}$, and the same argument shows that it fixes $\Gamma(x)\cup\Gamma(y)$ pointwise. Induction along paths and the connectedness of $\Gamma$ now give $G_{uv}=1$.

Consequently, $|G_u|=k$ and $k_2:=|\Gamma_2(u)|$ divides $k$. Moreover,
\[
2k=|G:F|\leq |\Aut(F)|=\varphi(m)\leq m-1.
\]
If $\Gamma$ had diameter $2$, it would be $\K_{m,m}$, with $k=m$, contradicting this inequality. Thus $d(\Gamma)\geq 3$. Since $\Gamma$ is amply regular, Theorem~\ref{Thm1.5.5} gives $k_2\geq k$, so $k_2=k$ and $\Gamma$ is a polygon or a Taylor graph.

If $k=2$, then $\Gamma=C_{2m}$ and $|G|=4m$, so $G=\Aut(\Gamma)$. The cyclic rotation subgroup of order $2m$ is normal in $G$ and hence contained in its Fitting subgroup $F$, contradicting $|F|=m$.

If $k\geq 3$, then $\Gamma$ is a bipartite Taylor graph, with intersection array $(k,k-1,1;1,k-1,k)$. By~\cite[Lemma~4.4(ii)]{CJS-2dt}, $\Gamma\cong\K_{k+1,k+1}-(k+1)\K_2$. Thus $m=k+1$, contradicting $2k\leq m-1=k$.

Altogether, $F$ cannot have two orbits on $V(\Gamma)$.
This completes the proof.
\qed

We are now ready to prove Theorem~\ref{2dt-lp-sf}.

\medskip
\noindent {\bf Proof of Theorem~\ref{2dt-lp-sf}.}
Assume that $\Gamma$ has a square-free number of vertices.
Let $F$ be the Fitting subgroup of $G$.
Since $\Gamma$ is $G$-locally-primitive, we have $\Gamma\ncong \K_{m[b]}$ with $m\geq 3$ and $b\geq 2$.
Therefore, by Lemma~\ref{2gt-sqf-soluble-1}, either $|F|$ is square-free, or $\Gamma\cong \K_{p,p}$, where $p$ is a prime. Therefore, to prove Theorem~\ref{2dt-lp-sf}, it remains to show that $|F|$ cannot be square-free.

Assume that $|F|$ is square-free.
Then it follows from Lemma~\ref{2dt-sqf-quot-1} that $F=C_G(F)$ is cyclic and regular on $V(\Gamma)$.
Thus $F$ is a normal regular subgroup of $G$, and hence we can identify $V(\Gamma)$ with $F$ such that $\Gamma\cong \Cay(F,S)$, where
$S=S^{-1}\subset F\setminus \{1\}$ and $F=\langle S\rangle$.
Let $u=1_F$. Since $F$ is normal and regular in $G$, $G_u\leq \Aut(F)$.

Since $S=S^{-1}$, for each $x\in S$, we have $\{x,x^{-1}\}\subseteq \Gamma(u)$.
As $G_u\leq \Aut(F)$, the set $\{x,x^{-1}\}$ is a block of the $G_u$-action on $\Gamma(u)$.
By the primitivity of $G_u$ on $\Gamma(u)=S$, the set $\{x,x^{-1}\}$ must be a trivial block.
Since the valency of $\Gamma$ is at least $3$, we have $|\{x,x^{-1}\}|=1$, and so $x=x^{-1}$.
Therefore each $x\in S$ has order $2$, and $S$ consists of involutions. As $F$ is cyclic, we have $F=\langle S\rangle\cong \mathbb{Z}_2$, which implies $|V(\Gamma)|=|F|=2$, a contradiction.

Therefore, $|F|$ cannot be square-free. This proves Theorem~\ref{2dt-lp-sf}.
\qed

\section*{Acknowledgements}
The authors would like to thank the anonymous referee for valuable comments and suggestions.
Wei Jin is  supported by NSFC (12271524,12331013),  NSF of Hunan (2026JJ50358) and Scientific Research Project of  the Education Department of Hunan Province (24A0142). Jack H. Koolen is supported by NSFC (12471335) and the Anhui Initiative in Quantum Information Technologies (AHY150000).

\begin{appendices}
\section{Two families of 2-arc-transitive graphs}\label{app-2at}
Let $q \geq 5$ be an odd prime power. In~\cite{JKL-2025-1}, the authors of the present paper constructed three families of amply regular graphs with diameter $d = 4$ and parameters $(4q+4, q, 0, \frac{q - 1}{2})$, arising from Paley graphs, Peisert graphs, and Paley digraphs, respectively. In this section, we prove that two of these families, namely those constructed from Paley graphs and Paley digraphs, are $2$-arc-transitive. Further examples of $2$-arc-transitive graphs arising from group divisible designs are discussed in~\cite{QKP-2016}.

\begin{prop}\label{prop:paley}
Let $\Delta$ be  the Paley graph $P(q)$, where $q$ is a prime power satisfying $q\equiv 1\pmod{4}$.
Let $\Sigma$ be the Taylor extension of $\Delta$, and let $\Gamma$ be the bipartite double of $\Sigma$. Then $\Gamma$ is $2$-arc-transitive.
\end{prop}
\proof
Let $\mathbb{F}_q$ be a finite field of order $q$.
Define the set of nonzero quadratic residues by
\[
\QR
=
\{x\in \mathbb{F}_q^\times \mid x=y^2 \text{ for some } y\in\mathbb{F}_q\},
\]
and the set of quadratic non-residues by
\[
\NQR
=
\mathbb{F}_q^\times\setminus \QR.
\]

Let $\Delta$ be the Paley graph $P(q)$, and let $\Sigma$ be the Taylor extension of $\Delta$. Then we may assume that
\[
V(\Sigma)
=
\{A_x,B_x \mid x\in\mathbb{F}_q\}
\cup
\{\infty_A,\infty_B\}.
\]
The edge set $E(\Sigma)$ is given by
\begin{align*}
E(\Sigma)
=
&\{\infty_A A_x \mid x\in\mathbb{F}_q\}
\cup
\{\infty_B B_x \mid x\in\mathbb{F}_q\}\\
&\cup
\{A_xA_y \mid x-y\in\QR\}
\cup
\{B_xB_y \mid x-y\in\QR\}\\
&\cup
\{A_xB_y \mid x-y\in\NQR\}.
\end{align*}
Let $\Gamma$ be the bipartite double of $\Sigma$. Then we may assume that
\[
V(\Gamma)
=
\{A_x^\sigma,B_x^\sigma \mid x\in\mathbb{F}_q,\ \sigma\in\{+,-\}\}
\cup
\{\infty_A^+,\infty_A^-,\infty_B^+,\infty_B^-\}.
\]
The edge set $E(\Gamma)$ is given by
\begin{align*}
E(\Gamma)
=
&\{\infty_A^-A_x^+ \mid x\in\mathbb{F}_q\}
\cup
\{\infty_A^+A_x^- \mid x\in\mathbb{F}_q\}\\
&\cup
\{\infty_B^-B_x^+ \mid x\in\mathbb{F}_q\}
\cup
\{\infty_B^+B_x^- \mid x\in\mathbb{F}_q\}\\
&\cup
\{A_x^+A_y^- \mid x-y\in\QR\}
\cup
\{B_x^+B_y^- \mid x-y\in\QR\}\\
&\cup
\{A_x^+B_y^- \mid x-y\in\NQR\}
\cup
\{A_x^-B_y^+ \mid x-y\in\NQR\}.
\end{align*}

Define a map $\varphi: V(\Gamma) \to V(\Gamma)$ by
\begin{align*}
\varphi(\infty_A^\sigma) &= A_0^{-\sigma}, & \varphi(A_0^\sigma) &= \infty_A^{-\sigma}, \\
\varphi(\infty_B^\sigma) &= B_0^{-\sigma}, & \varphi(B_0^\sigma) &= \infty_B^{-\sigma},
\end{align*}
where $\sigma\in\{+,-\}$ and $-\sigma$ denotes the opposite sign. For each $x\in\mathbb{F}_q^\times$, define
\begin{align*}
\varphi(A_x^\sigma)
=
\begin{cases}
A_{-1/x}^{-\sigma}, & x\in\QR,\\
B_{-1/x}^{-\sigma}, & x\in\NQR,
\end{cases}
\qquad
\varphi(B_x^\sigma)
=
\begin{cases}
B_{-1/x}^{-\sigma}, & x\in\QR,\\
A_{-1/x}^{-\sigma}, & x\in\NQR.
\end{cases}
\end{align*}
It is straightforward to verify that $\varphi$ is an automorphism of $\Gamma$.

For each $h\in\mathbb{F}_q$, define a map
\(
\varphi_h:V(\Gamma)\to V(\Gamma)
\)
by
\begin{align*}
\varphi_h(A_x^\sigma) &= A_{x+h}^\sigma, & \varphi_h(B_x^\sigma) &= B_{x+h}^\sigma,
\end{align*}
for all $x\in\mathbb{F}_q$ and $\sigma\in\{+,-\}$, and let
\[
\varphi_h(\infty_A^\sigma)=\infty_A^\sigma,
\qquad
\varphi_h(\infty_B^\sigma)=\infty_B^\sigma.
\]
It is straightforward to verify that $\varphi_h\in\mathrm{Aut}(\Gamma)$.

The automorphisms $\varphi$ and $\varphi_h$ together imply that $\Gamma$ is vertex-transitive.

Let $t$ be a generator of the multiplicative group $\mathbb{F}_q^\times$. 

For each integer $i$ with $1\le i\le \frac{q-1}{2}$, define a map
\(
\psi_i:V(\Gamma)\to V(\Gamma)
\)
by
\begin{align*}
\psi_i(A_x^\sigma) &= A_{xt^{2i}}^\sigma, & \psi_i(B_x^\sigma) &= B_{xt^{2i}}^\sigma,
\end{align*}
for all $x\in\mathbb{F}_q$ and $\sigma\in\{+,-\}$, and let
\[
\psi_i(\infty_A^\sigma)=\infty_A^\sigma,
\qquad
\psi_i(\infty_B^\sigma)=\infty_B^\sigma.
\]
Since $t^{2i}\in\QR$, the map $\psi_i$ preserves the defining adjacency relations of $\Gamma$. Hence $\psi_i\in\mathrm{Aut}(\Gamma)$.

The automorphisms $\varphi$, $\varphi_h$, and $\psi_i$ together imply that $\Gamma$ is arc-transitive.

Define a map
\(
\theta:V(\Gamma)\to V(\Gamma)
\)
by
\begin{align*}
\theta(A_x^+) &= A_{tx}^+, & \theta(B_x^+) &= B_{tx}^+, \\
\theta(A_x^-) &= B_{tx}^-, & \theta(B_x^-) &= A_{tx}^-,
\end{align*}
for all $x\in\mathbb{F}_q$ and $\sigma\in\{+,-\}$, and let
\[
\theta(\infty_A^-)=\infty_A^-,
\qquad
\theta(\infty_B^-)=\infty_B^-.
\]
\[
\theta(\infty_A^+)=\infty_B^+,
\qquad
\theta(\infty_B^+)=\infty_A^+.
\]
It is straightforward to verify that $\theta\in\mathrm{Aut}(\Gamma)$.

Combining the four types of automorphisms $\varphi$, $\varphi_h$, $\psi_i$, and $\theta$ constructed above, we conclude that $\Gamma$ is $2$-arc-transitive.
\qed

Let $q \geq 7$ be a prime power such that $q \equiv 3 \pmod{4}$. Then the finite field $\mathbb{F}_q$ contains no square root of $-1$. Consequently, for each pair of distinct elements $a, b \in \mathbb{F}_q$, exactly one of $a - b$ and $b - a$ is a square in $\mathbb{F}_q$.
The \emph{Paley digraph} $\vec{P}(q)$ is the directed graph with vertex set $\mathbb{F}_q$, where there is an arc from $a$ to $b$ (with $a \ne b$) if and only if $b - a \in 
\Box = \{ x \in \mathbb{F}_q^\times \mid x = y^2 \text{ for some } y \in \mathbb{F}_q \}$.

Following an argument similar to that in Proposition~\ref{prop:paley}, one can show that the graphs constructed from Paley digraphs in~\cite{JKL-2025-1} are $2$-arc-transitive. We therefore state the following proposition without proof.

\begin{prop}
Let $q \geq 7$ be a prime power such that $q \equiv 3 \pmod{4}$. 
Let $A_0 = I_q$, let $A_1$ be the adjacency matrix of $\vec{P}(q)$, and let $A_2 = A_1^T$ be its transpose. 
Let
\[
C = \begin{bmatrix}
0 & \mathbf{1} & 0 & 0 \\
0 & A_1 & A_2 & \mathbf{1}^T \\
\mathbf{1}^T & A_2 & A_1 & 0 \\
0 & 0 & \mathbf{1} & 0
\end{bmatrix}, \qquad
B = \begin{bmatrix}
0 & C \\
C^T & 0
\end{bmatrix}.
\]
Here $\mathbf{1}$ denotes the all-one row vector of length $q$.
Let $\Gamma$ be the undirected graph with adjacency matrix $B$.
Then $\Gamma$ is $2$-arc-transitive.
\end{prop}

\end{appendices}

\end{document}